\documentclass[12pt]{article}
\usepackage[russian]{babel}
\usepackage{amsthm,amsfonts,amssymb,amscd}

\begin{document}


\begin{center}
\Large \bf Birationally rigid \\
Fano complete intersections. II
\end{center}
\vspace{1cm}

\centerline{Aleksandr Pukhlikov}

\parshape=1
3cm 10cm \noindent {\small \quad \quad \quad
\quad\quad\quad\quad\quad\quad\quad {\bf }\newline We prove that
a generic (in the sense of Zariski topology) Fano complete
intersection $V$ of the type $(d_1,\dots,d_k)$ in ${\mathbb
P}^{M+k}$, where $d_1+\dots+d_k=M+k$, is birationally superrigid
if $M\geq 7$, $M\geq k+3$ and $\mathop{\rm max} \{d_i\}\geq 4$.
In particular, on the variety $V$ there is exactly one structure
of a Mori fibre space (or a rationally connected fibre space),
the groups of birational and biregular self-maps coincide,
$\mathop{\rm Bir} V= \mathop{\rm Aut} V$, and the variety $V$ is
non-rational. This fact covers a consider\-ably larger range of
complete intersections than the result of [J. reine angew. Math
{\bf 541} (2001), 55-79], which required the condition $M\geq
2k+1$.

Bibliography: 14 titles.} \vspace{1cm}

\begin{flushright}
{\it To the memory of Eckart Viehweg}
\end{flushright}
\vspace{1cm}

{\bf 1. Formulation of the main result.} Fix a set of integers
$d_i\geq 2$, $i=1,\dots,k$, where $k\geq 2$,
satisfying the following conditions:

\begin{itemize}

\item $\mathop{\rm max}\limits_{i=1,\dots,k} \{d_i\}\geq 4$,

\item $d_1+\dots+d_k\geq \mathop{\rm max} \{2k+3,k+7\}$.

\end{itemize}

Set $M=d_1+\dots+d_k-k$. Obviously, $M\geq 7$ and $M\geq k+3$.
By the symbol ${\mathbb P}$ we denote the complex projective space
${\mathbb P}^{M+k}$. The integers $d_i$ are assumed to be non-decreasing:
$d_i\leq d_{i+1}$.

Consider a smooth complete intersection
$$
V=F_1\cap\dots\cap F_k\subset{\mathbb P}
$$
of the type $(d_1,\dots,d_k)$, that is, $\mathop{\rm deg}F_i=d_i$,
$\mathop{\rm dim}V=M$. Obviously, $V$ is a smooth Fano variety, $\mathop{\rm Pic}V={\mathbb Z}H$, where $H$ is the class of a hyperplane section and $K_V=-H$, that is, $V$ is a Fano variety of index one. The main result of the present note is the following\vspace{0.1cm}

{\bf Theorem 1. } {\it A generic (in the sense of Zariski topology)
complete intersection $V$ of the type $(d_1,\dots,d_k)$ is birationally superrigid.}\vspace{0.1cm}

{\bf Corollary 1.} {\it There are no other structures of a Mori fiber space on $V$, apart from $V/\mathop{\rm pt}$: if $\chi\colon V\to V'$ is a birational map on the total space of a Mori fibre space $V'/S'$, then $\chi$ is a biregular isomorphism and $S'$ is a point. There exist no rational dominant maps $\alpha\colon V\dashrightarrow Y$ onto positive-dimensional varieties $Y$, the generic fibre of which is rationally connected. The groups of birational and biregular self-maps of the variety $V$ coincide:} $\mathop{\rm Bir}V=\mathop{\rm Aut}V$.\vspace{0.1cm}

The claim of Theorem 1 is proved in [1] for complete intersections, satisfying the condition $V\geq 2k+1$. The present paper essentially improves that result: birational superrigidity remains an open question only for complete intersections of quadrics and cubics $(2,\dots,2,3,\dots,3)$ and for one exceptional series $(2,\dots,2,4)$; for a fixed dimension $M$ there are $[M/2]+1$ such families. All the other Fano complete intersections of index one are covered by [1] and the present paper.\vspace{0.1cm}

{\it Birational superrigidity} in the formulation of Theorem 1 is understood in the sense of [1,2] and other papers of the author: for any mobile linear system $\Sigma\subset|nH|$, $n\geq 1$,
the equality
$$
c_{\rm virt}(\Sigma)=c(\Sigma, V)=n
$$
holds, where $c(\cdot)$ is the {\it threshold of canonical adjunction}, and
$c_{\rm virt}(\cdot)$ is the {\it virtual threshold of canonical adjunction}. Obtaining Corollary 1 from Theorem 1 is a well known elementary exercise, see, for instance [1-3].\vspace{0.1cm}

There are other definitions of birational (super)rigidity. Often (see, for instance, [4]) the first claim of Corollary 1 (the uniqueness of the structure of a Mori fibre space) is called birational rigidity, which together with the equality $\mathop{\rm
Bir}V=\mathop{\rm Aut}V$ is called birational superrigidity. For Fano varieties of index one with the Picard group  ${\mathbb Z}$ all the definitions that are currently in use are equivalent.

Now let us formulate the standard conjectures on birational geometry of Fano varieties. \vspace{0.1cm}

{\bf Conjecture 1.} {\it A smooth Fano complete intersection $X\subset
{\mathbb P}(a_0,\dots,a_{M+k})$ of codimension $k$ and index one in a weighted projective space is birationally superrigid for $M\geq 5$ and birationally rigid for  $M\geq 4$. A generic (in the sense of Zariski topology) complete intersection is birationally superrigid for} $M\geq 4$.\vspace{0.1cm}

Conjecture 1 is confirmed by all the known results on birational rigidity of Fano varieties, see [5-7] and the bibliography in [2]. Note that there is an example of a higher-dimensional Fano variety of index one with the Picard group ${\mathbb Z}$, which is birational to a Fano variety of high index and for that reason is not birationally rigid [8]. In that example the Fano variety has the numerical Chow group $A^2$ of codimension two cycles larger than ${\mathbb Z}$, and this difference between the example and the complete intersections seems to be crucial.\vspace{0.1cm}

The paper is organized in the following way. In Sec. 2 we give a precise meaning to the assumption that the complete intersection $V\subset{\mathbb
P}$ is generic and make the formulation of the main result more precise (Theorem 3). In Sec. 3 we prove birational superrigidity of complete intersections.\vspace{0.5cm}


{\bf 2. Generic Fano complete intersections.} For a point $o\in
{\mathbb P}$ fix the standard affine set ${\mathbb
C}^{M+k}\subset {\mathbb P}$ with a system of linear coordinates
$(z_1,\dots,z_{M+k})$, where $o=(0,\dots,0)$. Set
$$
f_i=q_{i,0}+q_{i,1}+\dots+q_{i,d_i},
$$
where the polynomials $q_{i,a}(z_*)$ are homogeneous of degree $a$. Obviously, $o\in
V(f_1,\dots,f_k)$, if and only if $q_{1,0}=\dots=q_{k,0}=0$.\vspace{0.1cm}

Let us introduce two sets of pairs of indices:
$$
J=\{(i,j)\,|\, 1\leq i\leq k,\,\,1\leq j\leq
d_i,\,\,(i,j)\neq(k,d_k)\}\subset{\mathbb Z}_+\times{\mathbb Z}_+
$$
and
$$
J^+=\{(i,j)\,|\,1\leq i\leq k,\,\,2\leq j\leq
d_i,\,\,(i,j)\neq(k,d_k)\}\subset{\mathbb Z}_+\times{\mathbb Z}_+.
$$
Recall [1] the following\vspace{0.1cm}

{\bf Definition 1.} The point $o\in V(f_1,\dots, f_k)$ is {\it
regular}, if the set of homogeneous polynomials
$$
Q=\{q_{i,j}, (i,j)\in J\}
$$
forms a regular sequence in ${\cal O}_{o,{\mathbb
C}^{M+k}}={\cal O}_{o,{\mathbb P}}$.\vspace{0.1cm}

A regular point is non-singular: the tangent space
$$
T_oV={\mathbb T}=\{q_{1,1}=\dots=q_{k,1}=0\}\subset{\mathbb C}^{M+k}
$$
is a linear subspace of codimension $k$. Therefore, the point
$o\in V(f)$ is regular, if and only if it is non-singular and the set of homogeneous polynomials
$$
Q^+=\{q_{i,j}|_{\mathbb T},\,(i,j)\in J^+\}
$$
forms a regular sequence in ${\cal O}_{o,{\mathbb T}}$. The latter condition is equivalent to the closed algebraic set
$$
\{q_{i,j}|_{\mathbb T}=0,\,(i,j)\in J^+\}
$$
being a finite set of points in ${\mathbb P}(\mathbb T)={\mathbb
P}^{M-1}$.\vspace{0.1cm}

{\bf Definition 2.} A regular point $o\in V(f_*)$ is {\it
correct in quadratic terms}, if none of the irreducible components of the closed set
\begin{equation}\label{sept2011_1}
\{q_{1,2}|_{\mathbb T}=\dots=q_{k,2}|_{\mathbb T}=0\}
\end{equation}
is contained in a linear subspace of codimension two in $\mathbb T$.\vspace{0.1cm}

In other words, the regular point $o$ is correct in quadratic terms, if for any irreducible component $W$ of the closed set (\ref{sept2011_1}) its linear span $<W>$ is either a hyperplane in $\mathbb T$, or the whole space $\mathbb T$. Since all polynomials
$q_{i,j}$ are homogeneous, Definition 2 can be understood in terms of the projective space ${\mathbb P}(\mathbb T)={\mathbb P}^{M-1}$.\vspace{0.1cm}

{\bf Theorem 2.} {\it For a generic (in terms of Zariski topology on the space of $k$-uples $(f_1,\dots,f_k)$) complete intersection $V(f_1,\dots,f_k)$ every point $o\in V$ is regular and correct in quadratic terms.}\vspace{0.1cm}

By Theorem 2, the main result of the present paper (Theorem 1) is implied by the following fact.\vspace{0.1cm}

{\bf Theorem 3.} {\it A complete intersection $V\subset{\mathbb P}$,
which is regular and correct in quadratic terms at every point is birationally superrigid.}\vspace{0.1cm}

{\bf Proof} of Theorem 3 is given below in Sec. 3.\vspace{0.1cm}

{\bf Proof of Theorem 2.} Fix a point $o\in{\mathbb
P}$, an affine subset ${\mathbb C}^{M+k}$ with coordinates $(z_1,\dots,z_{M+k})$, where $o=(0,\dots,0)$, and a set of linearly independent linear forms $q_{1,1}(z_*),\dots,q_{k,1}(z_*)$, so that $\mathbb T=\{q_{1,1}=\dots=q_{k,1}=0\}$ is a linear subspace of codimension $k$. Set
$$
{\cal L}_1=\{(q_{i,j}|_{\mathbb T}, (i,j)\in J^+)\}
$$
to be the space of tuples of $(M-1)$ homogeneous polynomials on ${\mathbb P}^{M-1}={\mathbb P}(\mathbb T)$,
$$
{\cal L}_2=\{(q_{i,2}|_{\mathbb T}, i=1,\dots,k)\}
$$
to be the space of tuples of $k$ quadratic polynomials on ${\mathbb P}^{M-1}$. Define the closed set $Y_1\subset{\cal L}_1$ by the condition that the set $(q_{i,j}|_{\mathbb T}, (i,j)\in J^+)$ does not form a regular sequence (that is, the set of zeros of those polynomials is of positive dimension). Define the closed set $Y_2\subset{\cal L}_2$ as the closure of the set $Y^o_2$ of such tuples of quadratic polynomials $(q_{i,2}|_{\mathbb T}, i=1,\dots,k)$, that
\begin{itemize}
\item  the set $Z(q_{*,2})$ of their common zeros is of codimension $k$ in ${\mathbb P}^{M-1}$,
\item there is an irreducible component $W$ of the set $Z(q_{*,2})$, the linear span $<W>$ of which is of codimension $\geq 2$ in ${\mathbb P}^{M-1}$.
\end{itemize}
It is easy to see (as in [1, Sec. 3.2]), that Theorem 2 follows immediately from\vspace{0.1cm}

{\bf Proposition 1.} {\it The following inequalities hold:} 
$$
\mathop{\rm codim}(Y_1\subset{\cal L}_1)\geq M+1,  \quad \mathop{\rm codim}(Y_2\subset{\cal L}_2)\geq M+1.
$$\vspace{0.1cm}

{\bf Proof.} Let us obtain the first estimate. Order the set of polynomials $q_{i,j}|_{\mathbb T}$ in some way, so that
$$
(q_{i,j}|_{\mathbb T},(i,j)\in J^+)=(p_1,\dots,p_{M-1}),
$$
we may assume that $\mathop{\rm deg}p_i\leq \mathop{\rm deg}p_{i+1}$, the polynomials $p_i$ are considered on ${\mathbb P}^{M-1}$.\vspace{0.1cm}

Let ${\cal L}_{1,l}=\{(p_1,\dots,p_l)\}$ be the space of truncated tuples, $Y_{1,l,a}$ the closure of the set $Y^o_{1,l,a}\subset{\cal L}_{1,l}$, defined by the following conditions:
\begin{itemize}
\item  the set $Z(p_1,\dots,p_{l-1})=\{p_1=\dots=p_{l-1}=0\}$ is of codimension $l-1$ in ${\mathbb P}^{M-1}$,
\item  there is an irreducible component $W$ of the set $Z(p_1,\dots,p_{l-1})$, on which $p_l$ vanishes, and the linear span $<W>$ is of codimension $a\leq l-1$ in
${\mathbb P}^{M-1}$.
\end{itemize}
It is easy to see that the first inequality of Proposition 1 follows from the estimates
\begin{equation}\label{sept2011_2}
\mathop{\rm codim}(Y_{1,l,a}\subset{\cal L}_{1,l})\geq M+1
\end{equation}
for any $l\leq M-1$ and $a\leq l-1$. Let us prove the inequality (\ref{sept2011_2}).  Consider first the case $a=l-1$. Here $W=<W>$ is a linear subspace of codimension $l-1$, and
$$
p_1|_W\equiv 0,\dots,p_l|_W\equiv 0.
$$
Since $\mathop{\rm deg}p_i\geq 2$ and the dimension of the corresponding Grassmanian is equal to $(l-1)(M-l+1)$, we obtain the estimate
$$
\mathop{\rm codim}(Y_{1,l-1,l-1}\subset{\cal L}_{1,l})\geq
l{M-l+2 \choose 2} -(l-1)(M-l+1)=
$$
$$
=\frac{M-l+1}{2}(lM-l^2+2).
$$
It is easy to check that the minimum of the latter expression is attained at $l=M-1$ and equal to $M+1$. This proves the inequality (\ref{sept2011_2}) for $a=l-1$.\vspace{0.1cm}

Consider the case $a\leq l-2$. Fixed a linear subspace $P\subset{\mathbb P}^{M-1}$ of codimension $a$ and set $Y^o_{1,l,a}(P)\subset Y^o_{1,l,a}$ to be the set of such tuples that there exists an irreducible component $W$ of the set $Z(p_1,\dots,p_l)$ such that $<W>=P$ and $p_l|_W\equiv 0$. Let $Y_{1,l,a}(P)\subset Y_{1,l,a}$ be the closure of the set $Y^o_{1,l,a}(P)$.\vspace{0.1cm}

Let us estimate the codimension of the set $Y_{1,l,a}(P)$ in ${\cal L}_{1,l}$. We use the method of the proof of Proposition 4 in [1]: there are $l-a-1$ polynomials
$$
p_{i_1},\dots,p_{i_{l-a-1}}
$$
such that $W$ is an associated subvariety of the set $(p_{i_1},\dots,p_{i_{l-a-1}})$ (see [1, Definition 3]), so that
$$
\mathop{\rm codim}(Y_{1,l,a}(P)\subset{\cal L}_{1,l})
\geq(a+1)(2(M-a)-1),
$$
and thus
$$
\mathop{\rm codim}(Y_{1,l,a}\subset{\cal
L}_{1,l})\geq(a+1)(2(M-a)-1)-a(M-a)=
$$
$$
=a(M-1-a)+2(M-a)-1.
$$
The minimum of this function, quadratic in $a$, on the interval $[0,l-2]$ is attained at one of its endpoints. For $a=0$ we get $2M-1$. For $a=l-2$ we get
$$
l(M-l+1)+1.
$$
The minimum of the latter function on the interval $[2,M-1]$ is attained at its endpoints and equal to $2M-1$. This completes the proof of the first inequality.\vspace{0.1cm}

Let us prove the second inequality. We follow the same scheme as above. Let us define the set $Y_{2,l,a}$ as the closure of the set $Y^o_{2,l,a}$ in the space of truncated tuples ${\cal L}_{2,l}=\{(p_1,\dots,p_l)\}$, $l\leq k$, which is determined by the following condition: there exists an irreducible component $W$ of the set $Z(p_1,\dots,p_l)$, the linear span of which is of codimension $a\geq 2$. Obviously, it is sufficient to prove the inequality
$$
\mathop{\rm codim}(Y_{2,l,a}\subset{\cal L}_{2,l})\geq M+1
$$
for each $a=2,\dots,l$. Consider first the case $a=l$. Here $W=<W>$ is a linear space and we get the estimate
$$
\mathop{\rm codim}(Y_{2,l,l}\subset {\cal L}_{2,l})=l\frac{(M-l)(M-l+1)}{2}-l(M-l)=
$$
$$
=l(M-l)\frac{M-l-1}{2},
$$
which is by far stronger than we need.\vspace{0.1cm}

Consider the case $a\leq l-1$. Fix a linear subspace $P\subset{\mathbb P}^{M-1}$ of codimension $a$ and construct the set $Y^o_{2,l,a}(P)$ of such tuples $(p_1,\dots,p_l)$, that the closed set $Z(p_1,\dots,p_l)$ has an irreducible component $W$ such that $<W>=P$. Since $\mathop{\rm codim}W=l$, among the polynomials $p_1,\dots,p_l$ we can find
$(l-a)$ quadrics $p_{i_1},\dots,p_{i_{l-a}}$, such that $p_{i_1}|_P,\dots,p_{i_{l-a}}|_P$ form a good sequence and $W$ is its associated subvariety (see the proof of Proposition 4 in [1]). This implies the estimate
$$
\mathop{\rm codim}(Y_{2,l,a}(P)\subset{\cal L}_{2,l})\geq
a(2(M-a-1)+1),
$$
so that
$$
\mathop{\rm codim}(Y_{2,l,a}\subset{\cal L}_{2,l}\geq a(2(M-a)-1)-a(M-a)=a(M-a-1).
$$
It is easy to check that the minimum of the right hand side on the interval $[2,l-1]$ is attained at $a=2$ and equal to $2(M-3)\geq M+1$ for $M\geq 7$. This proves the second inequality of Proposition 1. Proof of Theorem 2 is complete.\vspace{0.5cm}


{\bf 3. Proof of birational superrigidity.} Starting from this moment, we fix a complete intersection $V(f_1,\dots,f_k)\subset{\mathbb P}^{M+k}$, satisfying the regularity condition and the conditions of being correct in the quadratic terms at every point. Assume that the variety $V$ is not birationally superrigid. This implies in an easy way (see [1-3]), that on $V$ there exists a mobile linear system $\Sigma\subset|n H|$ with a {\it maximal singularity}, that is, there are: a (non-singular) projective model $V^{\sharp}$, a birational morphism $\varphi\colon V^{\sharp}\to V$ and an exceptional divisor $E^{\sharp}\subset V^{\sharp}$ such that the {\it Noether-Fano inequality}
$$
\mathop{\rm ord}\nolimits_{E^{\sharp}}\varphi^{*}\Sigma>na(E^{\sharp},V)
$$
holds, that is, $E^{\sharp}$ realizes a non-canonical singularity of the pair $(V,\frac{1}{n}\Sigma)$.\vspace{0.1cm}

The irreducible subvariety $B=\varphi(E^{\sharp})\subset V$ is the {\it centre} of the maximal singularity $E^{\sharp}$. It is easy to check that $\mathop{\rm mult}_B\Sigma> n$.\vspace{0.1cm}

There are three options for the codimension of the subvariety $B$:

(i) $\mathop{\rm codim}B=2$,

(ii) $\mathop{\rm codim}B=3$,

(iii) $\mathop{\rm codim}B\geq 4$,

\noindent and we must show that none of them takes place. If $\mathop{\rm codim}B=2$, then word for word the same arguments as in [1, Sec. 1.1] give a contradiction. Those arguments work for any variety $V$ (without any conditions of general position), satisfying the equality $A^2V={\mathbb Z}H^2$, where $H$ is the positive generator of the Picard group, and $A^2$ is the numerical Chow group of codimension two cycles. (The case $\mathop{\rm codim}B=2$ can also be excluded by the argument used below in codimension three.)

If $\mathop{\rm codim}B=3$, then we come to a contradiction with the following fact.\vspace{0.1cm}

{\bf Proposition 2.}  {\it For any subvariety $W\subset V$ of dimension $\geq k$ the following inequality holds:}
$$
\mathop{\rm mult}\nolimits_W\Sigma\leq n.
$$

{\bf Proof.} This is Proposition 3.6 in [9]. (Recall that $\mathop{\rm dim}V=M\geq k+3$ by assumption, so that if $\mathop{\rm codim}B=3$, then $\mathop{\rm dim}B\geq k$ and [9, Proposition 3.6] applies.)\vspace{0.1cm}

Therefore, we may assume that the third case takes place: $\mathop{\rm codim}B\geq 4$. Moreover, we may assume that the codimension of the subvariety $B$ is minimal among all centres of maximal singularities of the linear system $\Sigma$, in particular, $B$ is not contained in the centre $B'$ of another maximal singularity $E'$, if they do exist. Let $o\in B$ be a point of general position, $\lambda\colon\widetilde{V}\to V$ its blow up, $E=\lambda^{-1}(o)\cong{\mathbb P}^{M-1}$ the exceptional divisor. By the symbol $\widetilde{\Sigma}$ we denote the strict transform of the mobile system $\Sigma$ on $\widetilde{V}$. Let
$$
Z=(D_1\circ D_2)
$$
be the {\it self-intersection} of the system $\Sigma$, that is, the scheme-theoretic intersection of two generic divisors $D_1,D_2\in\Sigma$. By the symbol $\widetilde{Z}$ we denote the strict transform of the effective cycle $Z$ on $\widetilde{V}$. The following fact is true.\vspace{0.1cm}

{\bf Proposition 3.} {\it There exist a linear subspace $P\subset E$ of codimension two, satisfying the inequality
$$
\mathop{\rm mult}\nolimits_oZ+\mathop{\rm mult}\nolimits_P\widetilde{Z}>8n^2.
$$
If $\mathop{\rm mult}_oZ\leq 8n^2$, then the linear subspace $P$ is uniquely determined by the linear system $\Sigma$}.\vspace{0.1cm}

{\bf Proof.} Let $X\ni o$ be a generic germ of a smooth subvariety of dimension $\mathop{\rm codim}B\geq 4$. By the assumption on minimality of $\mathop{\rm codim}B$, the pair $(X,\frac{1}{n}\Sigma_X)$ is canonical outside the point $o$, but not canonical at that point, where $\Sigma_X=\Sigma|_X$ is the restriction of the mobile system $\Sigma$ onto $X$. Now, applying [10, Proposition 4.1], we obtain the required claim.\vspace{0.1cm}

Let us consider the tangent hyperplanes $T_oF_i$, $i=1,\dots,k$ (in the coordinates $z_*$ they are given by the equations $q_{i,1}=0$). Let $T_i=V\cap T_oF_i$ be the corresponding tangent divisors, $\mathop{\rm mult}_oT_i=2$ (the equality holds because of the regularity condition), and let $\widetilde{T_i}\subset\widetilde{V}$ be their strict transforms on $\widetilde{V}$. The equation $q_{i,2}|_{\mathbb T}=0$ defines the projectivised tangent cone $T^E_i=(\widetilde{T_i}\circ E)$, since $E={\mathbb P}({\mathbb T})\cong{\mathbb P}^{M-1}$. Since none of the irreducible components of the effective cycle $(T^E_1\circ\dots\circ T^E_k)$ of codimension $k$ is contained in a linear space of codimension two in $E$, for a generic hyperplane $\Lambda\supset P$ in $E$ we get
\begin{equation}\label{sept2011_3}
\mathop{\rm codim}\nolimits_E(T^E_1\cap\dots\cap T^E_k\cap\Lambda)=k+1.
\end{equation}

Let $L\in |H|$, $L\ni o$ be a generic hyperplane section, such that $\widetilde{L}\cap E=\Lambda$.  By genericity, none of the components of the cycle $Z$ is contained in $L$, so that $Z_L=(Z\circ L)$ is an effective cycle of codimension three on $V$, satisfying the estimate
$$
\mathop{\rm mult}\nolimits_oZ_L\geq\mathop{\rm mult}\nolimits_oZ+\mathop{\rm mult}\nolimits_P\widetilde{Z}>8n^2,
$$
and $\mathop{\rm deg}Z_L=\mathop{\rm deg}Z=dn^2$. Therefore, there exists an irreducible subvariety $C\subset V$ of codimension three (an irreducible component of the cycle $Z_L$), such that
$$
\frac{\mathop{\rm mult}_oC}{\mathop{\rm deg}C}>\frac{8}{d}.
$$
Obviously, the support of the projectivised tangent cone $C_E=(\widetilde{C}\circ E)$ is contained in the hyperplane $\Lambda$. Let $\Delta=\{(\lambda_1q_{1,2}+\dots+\lambda_kq_{k,2})|_{\mathbb T}=0\}$ be the linear system of quadrics on $E$, spanned by the quadrics $T^E_1,\dots,T^E_k$. 
By the equality (\ref{sept2011_3}), for $(k-2)$ generic divisors in this system, for simplicity of notations let them be just $T^E_1,\dots,T^E_{k-2}$, such that
$$
\mathop{\rm codim}\nolimits_E(C_E\cap T^E_1\cap\dots\cap T^E_{k-2})=k+1.
$$
Therefore, for the codimension in a neighborhood of the point $o$ we also have
$$
\mathop{\rm codim}\nolimits_o(C\cap T_1\cap\dots\cap T_{k-2})=k+1.
$$
Now let us construct a sequence of irreducible subvarieties $R_0,R_1,\dots,R_{k-2}$ with the following properties:\vspace{0.1cm}

(i) $R_0=C$, $R_i\ni o$ for all $i=1,\dots,k-2$;\vspace{0.1cm}

(ii) $R_{i+1}$ is an irreducible component of the effective cycle $(R_i\circ T_{i+1})$, for which the ratio
$$
\frac{\mathop{\rm mult}_o}{\mathop{\rm deg}}
$$
takes the maximal value.\vspace{0.1cm}

Obviously, if $R_i\ni o$ and $R_{i+1}$ satisfies the property (ii), then $\mathop{\rm mult}_oR_{i+1}>0$, so that $R_{i+1}\ni o$ as well. Moreover, $\mathop{\rm codim}_o(C\cap T_1\cap\dots\cap T_{i+1})=i+4$ and $\mathop{\rm codim}R_i=i+3$, so that $R_1\not\subset T_{i+1}$ and the construction described above is possible. It is easy to see that for all $i=1,\dots,k-2$
$$
\frac{\mathop{\rm mult}_o}{\mathop{\rm deg}}R_i\geq 2^i\cdot\frac{\mathop{\rm mult}_o}{\mathop{\rm deg}}C.
$$
Set $R=R_{k-2}$. This is an irreducible subvariety of codimension $k+1$, $o\in R$ and the inequality
$$
\frac{\mathop{\rm mult}_o}{\mathop{\rm deg}}R>\frac{2^{k+1}}{d}
$$
holds. Therefore, by [1, Corollary 1], we obtain the estimate
$$
\frac{2^{k+1}}{d}<\frac{\mathop{\rm mult}_o}{\mathop{\rm deg}}R\leq\frac{3d_k}{2d_k-2}\cdot 2^b\cdot\left(\prod_{d_i\geq 3}d_i\right)^{-1},
$$
where $b=k-{\sharp}\{i\,|\, 1\leq i\leq k,d_i=2\}$. The proof, given in [1, Sec. 2.2-2.3] does not make use of the condition $M\geq 2k+1$ and so works for any complete intersection, satisfying the regularity condition. Since
$$
d=2^{k-b}\left(\prod_{d_i\geq 3}d_i\right),
$$
we obtain the inequality
$$
1<\frac{3d_k}{4d_k-4},
$$
which is not true for $d_k\geq 4$. This completes the proof of Theorem 3, and so that of Theorem 1. Q.E.D.\vspace{0.1cm}

{\bf Remark 1.} The key difference of the arguments of the present paper from the proof given in [1] is that the combination of the $4n^2$-inequality with the Lefschetz theorem ($A^iV={\mathbb Z}H^i$ for $i>[M/2]$) is replaced by the $8n^2$-inequality. This makes it possible to avoid the assumption that $k<[M/2]$, which is needed to apply the Lefschetz theorem. The $4n^2$-inequality [2,3,11] goes back to the classical paper of V.A.Iskovskikh and Yu.I.Manin on the three-dimensional quartic [12]. The $8n^2$-inequality was known since 2000 (see [4]), however, its proof in [4,13] generated some doubts and, as it turned out, contained an essential gap indeed, which was corrected only very recently [10,14].\vspace{0.5cm}


{\bf References}\vspace{0.1cm}

{\small
\noindent 1. Pukhlikov A.V., Birationally rigid Fano complete
intersections, Crelle J. f\" ur die reine und angew. Math. {\bf
541} (2001), 55-79. \vspace{0.1cm}

\noindent 2. Pukhlikov A.V., Birationally rigid varieties. I.
Fano varieties. Russian Math. Surveys {\bf 62} (2007), No. 5,
857-942. \vspace{0.1cm}

\noindent 3. Pukhlikov A.V., Birational automorphisms of Fano
hypersurfaces, Invent. Math. {\bf 134} (1998), No. 2, 401-426.
\vspace{0.1cm}

\noindent 4. Iskovskikh V.A., Birational rigidity of Fano
hypersurfaces in the framework of Mori theory. Russian Math.
Surveys. {\bf 56} (2001), No. 2, 207-291.\vspace{0.1cm}

\noindent 5. Pukhlikov A.V., Birationally rigid iterated Fano
double covers. Izvestiya: Mathematics. {\bf 67} (2003), no. 3,
555-596. \vspace{0.1cm}

\noindent 6. Cheltsov I.A. Birationally rigid Fano varieties.
Russian Math. Surveys. {\bf 60} (2005), No. 5, 875-965.\vspace{0.1cm}

\noindent 7. Pukhlikov A.V., Birational geometry of algebraic
varieties with a pencil of Fano cyclic covers. Pure and Appl.
Math. Quart. {\bf 5} (2009), No. 2, 641-700.\vspace{0.1cm}

\noindent 8. Castravet A.-M., Examples of Fano varieties of index one that are not birationally rigid, Proc. Amer. Math. Soc. {\bf 135} (2007), No. 12, 3783-3788. \vspace{0.1cm}

\noindent 9. Pukhlikov A.V., Birational geometry of algebraic
varieties with a pencil of Fano complete intersections,
Manuscripta Mathematica {\bf 121} (2006), 491-526. \vspace{0.1cm}

\noindent 10. Pukhlikov A.V., Birational geometry of Fano double spaces of index two, Izvestiya: Mathematics {\bf 74} (2010), No. 5, 925-991.\vspace{0.1cm}

\noindent 11. Pukhlikov A.V., Essentials of the method of maximal
singularities, in ``Explicit Birational Geometry of Threefolds'',
London Mathematical Society Lecture Note Series {\bf 281} (2000),
Cambridge University Press, 73-100. \vspace{0.1cm}

\noindent 12. Iskovskikh V.A. and Manin Yu.I., Three-dimensional
quartics and counterexamples to the L\" uroth problem, Math. USSR
Sb. {\bf 86} (1971), No. 1, 140-166.\vspace{0.1cm}

\noindent 13. Cheltsov I., Double cubics and double quartics,
Math. Z. {\bf 253} (2006), No. 1, 75-86. \vspace{0.1cm}

\noindent 14. Pukhlikov A.V., On the self-intersection of a movable linear system,
Journ. of Math. Sci.  {\bf 164} (2010), No. 1, 119-133. \vspace{0.1cm}

}

\begin{flushleft}
{\it pukh@liv.ac.uk,\, pukh@mi.ras.ru}
\end{flushleft}

\end{document}